\documentclass[letterpaper, 10 pt, conference]{ieeeconf}  

\IEEEoverridecommandlockouts                              
\usepackage{graphicx} 
\usepackage{amsmath} 
\usepackage{amssymb}  
\usepackage{amstext,enumerate}

\newcommand{\zono}[1]{\langle #1 \rangle}

\usepackage{latexsym, color,booktabs}
\newtheorem{remark}{Remark}
\newtheorem{assumption}{Assumption}

\newtheorem{definition}{Definition}
\newtheorem{lemma}{Lemma}
\newtheorem{proposition}{Proposition}
\newtheorem{corollary}{Corollary}
\newtheorem{theorem}{Theorem}

\title{\LARGE \bf
Reachable set-based dynamic quantization for \\ the remote state estimation of linear systems
}
\author{Yaodong Li, Michelle S. Chong
\thanks{This work was done when Y. Li was a student at the Department of Mechanical Engineering, Eindhoven University of Technology. M. Chong is with the same department.  
        Emails:{\tt\small lyd1106nl@gmail.com, m.s.t.chong@tue.nl} }
\thanks{The authors thank K.J.A. Scheres for the initial discussions.}
}

\begin{document}
\maketitle
\thispagestyle{empty}
\pagestyle{empty}

\begin{abstract}
We employ reachability analysis in designing dynamic quantization schemes for the remote state estimation of linear systems over a finite date rate communication channel. The quantization region is dynamically updated at each transmission instant, with an approximated reachable set of the linear system. We propose a set-based method using zonotopes and compare it to a norm-based method in dynamically updating the quantization region. For both methods, we guarantee that the quantization error is bounded and consequently, the remote state reconstruction error is also bounded. To the best of our knowledge, the set-based method using zonotopes has no precedent in the literature and admits a larger class of linear systems and communication channels, where the set-based method allows for a longer inter-transmission time and lower bit rate. Finally, we corroborate our theoretical guarantees with a numerical example. 
\end{abstract}

\section{Introduction} \label{sec:intro}
 Cyber-physical systems integrate multiple agents and their sensing and actuation devices over a communication channel. Therefore, the remote estimator is responsible for reconstructing the state information based on the sensor's data transmitted over a communication channel. The remote state estimator is commonly used in control systems including autonomous vehicles, smart grids, and industrial automation \cite{hespanha2007survey}. For systems requiring signal communication via a bandwidth-limited network, the analog signal must be converted into discrete-valued digital symbols before being transmitted. This operation inevitably causes an error, called the quantization error.  Furthermore, the total amount of information that may be transmitted per unit of time is often limited due to bandwidth constraints of the digital channels, which further degrades the precision of the information that is exchanged over the network. There are results that show that the quantization error does not behave like, e.g., white noise \cite{delchamps1990stabilizing}. Hence, the existing observer designs such as {Kalman filters} may not be capable to handle the ``unknown noise'' caused by the quantization process. Therefore, novel approaches for remote state estimation over communication channels with bandwidth constraints are needed.

First, to minimize or even eliminate the effect of the quantization error,  many studies of quantizer design have been done in the past decades.  A quantizer is a mathematical mapping from a continuous region called the quantization region to a finite discrete set of indexes, which we call quantization levels. A quantizer with fixed parameters is called a static (or memory-less) quantizer.  In this case, a low bit rate has a strong negative effect on the resolution of quantization. Another strategy for improving the resolution under a fixed bandwidth constraint is \emph{dynamic quantization} \cite{brockett2000quantized}, where the quantization parameters are dynamically adjusted based on the received data and knowledge of the plant dynamics.  In \cite{liberzon2007input,sharon2011input}, the robustness of dynamic quantization against external disturbances is studied based on the notion of input-to-state stability \cite{sontag2008input}. The surveys \cite{nair2007feedback,jiang2013quantized} give recent developments of dynamic quantization for linear and nonlinear systems, respectively.

So far, stabilization problems with controller and quantizer designs are the main focus of the aforementioned papers, where the origin is often assumed to be a stable fixed point under the designed closed-loop controllers. On the other hand, the input and/or control laws are unknown to the remote state estimator. In this case, a conservative prediction for all the reachable states based on a constrained input is more reliable for adjusting the quantization parameters and resolution.   In \cite{combastel2005state,scott2016constrained,su2017model}, a set-based observer without the consideration of quantization effects propagates the set of all possible states  with bounded inputs using reachability analysis based on zonotopes \cite{girard2005reachability}. We will use this result to propagate a conservative approximation of the quantization region for dynamic quantization.
In this case, we choose the uniform quantizer \cite{de2004stabilizability} so that the resolution is evenly distributed across all components of the state. 

In this paper, we design dynamic quantization schemes by performing reachability analysis of continuous-time linear time-invariant (LTI) dynamical systems with external disturbances, for remote state estimation. We over-approximate the terminal reachable set at each transmission instant to update the quantization region and further improve the resolution of the quantization.  Polytope-based over-approximation of the reachable set can provide better results but the dramatic increase of vertices and surfaces makes it difficult to propagate. To that end, we were inspired by \cite{girard2005reachability} which uses \emph{zonotopes} to over-approximate the terminal reachable set of a LTI system which have desirable properties. The propagation of zonotopes relies on the \emph{centroid} and \emph{generators}, which are easy to compute and store. Although the number of generators increases linearly with time, in the dynamic quantization scheme, we only over-approximate
the terminal set in one inter-transmission interval, by which the issue is avoided. We then compare our results to a norm-based method, which was inspired by \cite{de2004stabilizability} that  uses the Lipschitz condition to upper-bound the terminal state. However, this only provides the norm bound instead of the more precise bound afforded by considering a component-wise bound of each individual state provided by the zonotopic method.

Similar to \cite{xu2005estimation}, we employ a pre-estimator before transmission to avoid remote estimation based on the outputs. We assume the input signal is known to the pre-estimator before transmission, but the remote state estimator has no access to either the inputs or the control laws. Under the aforementioned setup and its corresponding properties, several objectives can be achieved:
\begin{enumerate}
\item To the best of our knowledge, the first set-based dynamic quantization scheme using zonotopes.
\item Comparison between our set-based and a norm-based dynamic quantization scheme.
\item No overflow occurs for both dynamic quantization schemes.
\item Conditions for the transmission bit rate and inter-transmission interval are given to ensure the boundedness of the quantization error.
\item An upper bound is given for the state reconstruction error with respect to bounds on the quantization error, the input, and the disturbance.
\end{enumerate}

The paper is organized as follows. The next section presents the notation and other necessary preliminaries, definitions, and properties. In Section \ref{sec:prob}, we introduce the remote state estimation setup as well as its individual components separately and proceed to state the overall objective and approach taken in this paper. In Sections \ref{sec:dynamic_set} and \ref{sec:dynamic_norm}, we present two dynamic quantization schemes that achieve no overflow based on different propagation techniques. A comparison of the two methods is done in Section \ref{sec:compare}. Lastly, in Section \ref{sec:sim}, numerical simulations and comparisons are presented to support our results. We then conclude the paper with Section \ref{sec:conclude}.

\section{Preliminaries} \label{sec:prelim}

\subsection{Notations}\label{sec:notation}
\begin{itemize}
    \item Let $\mathbb{R}:=(-\infty,\infty)$, $\mathbb{R}_{\geq 0}:=[0,\infty)$, $\mathbb{Z}_{+}:=\{0,1,\dots\}$, $\mathbb{Z}_{>0}:=\{1,2,\dots\}$  and $\mathcal{I}_n:= \{1,2,...,n\}$.
    \item $\lambda_i(A)$ denotes the $i$-{th} eigenvalue of matrix $A$, $\lambda_{max}(A),\ \lambda_{min}(A)$ denote the maximum and the minimum eigenvalues, respectively.
    \item $I_n$ denotes the identity matrix of size $n \times n$ and $\textbf{0}_n$ denotes the $n$-dimensional vector $[0,0,...,0]^T$, and $\textbf{1}_{n\times m}$ denotes the $m\times n$ matrix with all elements as $1$.
    \item $\sigma_{n,i}$ is the n-dimensional vector with 1 as its $i$-th element and 0 otherwise. We write $\sigma_i$ when its dimension is clear from context.
    \item An $n$-dimensional hypercube with center $c\in\mathbb{R}^{n}$ and radius $l\in\mathbb{R}_{\geq 0}$ is denoted $\mathcal{B}(c,l)$.
    \item The infinity norm of a vector $x\in \mathbb{R}^n$ is denoted by $|x| : = \underset{{i\in \mathcal{I}_n}}\max{|x_i|}$. For a matrix $A \in \mathbb{R}^{n\times n}$, $|A| := \underset{{j\in \mathcal{I}_n}}\max{\sum_{i=1}^n{|a_{ij}|}}$, where $a_{ij}$ is the $(i,j)$-component of the matrix $A$.
    
    \item The Minkowski sum is denoted by $\oplus$. 
    \item $\lfloor b \rfloor$ denotes the floor function that returns the greatest integer that is less than or equal to $b \in \mathbb{R}$. 
    \item A continuous function $\gamma: \mathbb{R}_+\to \mathbb{R}_+$ is a class-$\mathcal{K}$ ($\mathcal{L}$) function, if it is strictly increasing (decreasing) and $\gamma(0) = 0$.  
    \item A continuous function $\alpha: \mathbb{R}_+\times \mathbb{R}_+ \to \mathbb{R}_+$ is a class-$\mathcal{KL}$ function, if $\alpha(\cdot,s)$ is a class-$\mathcal{K}$ for all $s\geq 0$, $\alpha(r,\cdot)$ is non-increasing and $\alpha(r,s) \to 0$ as $s \to \infty$ for all $r \geq 0$.
\end{itemize}

\subsection{Set representations} \label{sec:set}
A zonotope $\mathcal{Z}\subset \mathbb{R}^n$ is a set that satisfies
\begin{equation}
    \label{def zono}
    \mathcal{Z} := \left\{ x\in \mathbb{R}^n \bigg\vert x = c + \sum_{i = 1}^p \epsilon_i g_i,\ \forall \epsilon_i \in [-1,1] \right\},
\end{equation}
where $c\in \mathbb{R}^n$ is the geometric center of the zonotope and the line segments $g_i \in \mathbb{R}^n$ are called the generators of  the zonotope. We denote $\mathcal{Z} = \zono{c,G}$, where $G:=(g_1,...,g_p)\in\mathbb{R}^{n\times p}$.
Zonotopes possess the following properties:
\begin{enumerate}
    \item \underline{Minkowski sum}: The Minkowski sum of two zonotopes $\mathcal{Z}_a:=\zono{c_a,G_a}$ and $\mathcal{Z}_b:=\zono{c_b,G_b}$ remains a zonotope $\mathcal{Z}_c$ with $\mathcal{Z}_c := \mathcal{Z}_a \oplus \mathcal{Z}_b = \zono{c_a+c_b, (G_a,G_b)}$.
\item \underline{Linear transformation}: A matrix $K\in\mathbb{R}^{n\times n}$ multiplied with a zonotope $\mathcal{Z}:=\zono{c,G}$ results in a linearly transformed zonotope $K\mathcal{Z}=\zono{Kc,KG}$.
\item \underline{Containment of a zonotope by a hyperrectangle}: For any $x \in \mathcal{Z} = \zono{\mathbf{0},G}$, where $G:=(g_1,g_2,\dots,g_p)$, its $j$-th component $x_j$ is bounded in absolute value by $\sum_{i=1}^p|\sigma_j^T g_i|$, where $\sigma_j$ is a basis vector defined in Section \ref{sec:notation}. Therefore, a zonotope $\mathcal{Z}$ can be contained by a \emph{hyperrectangle} denoted by $\mathcal{H}(\mathbf{0},H)$, with $H:=(h_1,h_2,\dots,h_n)^T$, where $h_i:=\sum_{l=1}^{p}|\sigma_i^Tg_l|$ for $i\in\mathcal{I}_n$. 
\item \underline{Hyperrectangle in zonotopic form}: A hyperrectangle $\mathcal{H}(c,l)$ with $c\in\mathbb{R}^{n}$, $l=(l_1,l_2,\dots,l_n)^T$ can equivalently be written in zonotopic form as $\zono{c,(l_1\sigma_1,l_2\sigma_2,\dots,l_n\sigma_n)}=\zono{c,(\sigma_1^Tl\sigma_1,\sigma_2^Tl\sigma_2,\dots,\sigma_n^Tl\sigma_n)}$, where we obtain the last zonotopic form because $l_i=\sigma_i^T l$.
\end{enumerate}

\subsection{Terminal reachable set and its over-approximation}

The terminal reachable set of a dynamical system is defined as follows.

\begin{definition}
Consider a general dynamical system $\dot x = f(x,u)$, for a finite time interval $[t_0,t_1]$, the terminal reachable set $\mathcal{R}_{[t_0,t_1]}(\mathcal{X},\mathcal{U})$ is defined as the set of all states that are reachable at time $t_1$ with $x(t_0) \in \mathcal{X}$, with $u(t) \in \mathcal{U}$ for all $t\in [t_0, t_1]$, i.e.,
\begin{align}\label{def-reach}
    &\mathcal{R}_{[t_0,t_1]}\!(\mathcal{X},\mathcal{U})\!\notag \\
    &\;:=\!\{x(t_1)\!:\!\forall x(t_0)\!\in\!\mathcal{X}\!, \!u(t)\!\in\!\mathcal{U}\!, \dot{x}\!=\!f(x,u)\!, \forall t\!\in\![t_0,t_1]\!\}. \notag
\end{align}
\end{definition}

Next, we generalise Lemma 1 in \cite{girard2005reachability} to arbitrary finite time intervals $[t_0,t_1]$ for $(t_0,t_1)\in\mathbb{R}^2_{\geq 0}$ which we state below and will form a crucial step in the reachable set based approach we will take for our dynamic quantization schemes in this paper. 

\begin{lemma} \label{2009}
Consider a LTI system in the finite time interval $t\in[t_0,t_1]$, where $t_1-t_0 \leq \tau$ and $(t_0,t_1)\in \mathbb{R}^{2}_{\geq 0}$, given by
\begin{equation}
    \dot x(t) = Ax(t) + Bu(t), 
\end{equation}
with $x(t_0) \in \mathcal{X}, Bu(t) \in \mathcal{U}, \forall t\in [t_0, t_1]$. The sets $\mathcal{X}$ and $\mathcal{U}$ are zonotopes. Its terminal reachable set $\mathcal{{R}}_{[t_0,t_1]}(\mathcal{X},\mathcal{U})$ satisfies 
\begin{equation}
\label{reach}
    \mathcal{R}_{[t_0,t_1]}(\mathcal{X},\mathcal{U}) \subseteq e^{{\tau}A}\mathcal{X} \oplus \mathcal{B}(\textbf{0}_n,\beta(\tau,\mu)),
\end{equation}
where $\mu := \sup_{Bu\in\mathcal{U}}|Bu|$, and
\begin{equation}\label{betaT}
    \beta(\tau,\mu) : = |A|^{-1}e^{\tau|A|} \mu.
\end{equation}
\end{lemma}
Lemma \ref{2009} is a generalization of \cite[Lemma 1]{girard2005reachability} to an arbitrary finite time interval, which can be achieved through standard calculations based on the development in \cite{girard2005reachability}. Hence, its proof is omitted.

\section{Problem setup and approach} \label{sec:prob}
\subsection{Setup and plant}
\begin{figure}[htbp]
\centerline{\includegraphics[scale=0.28]{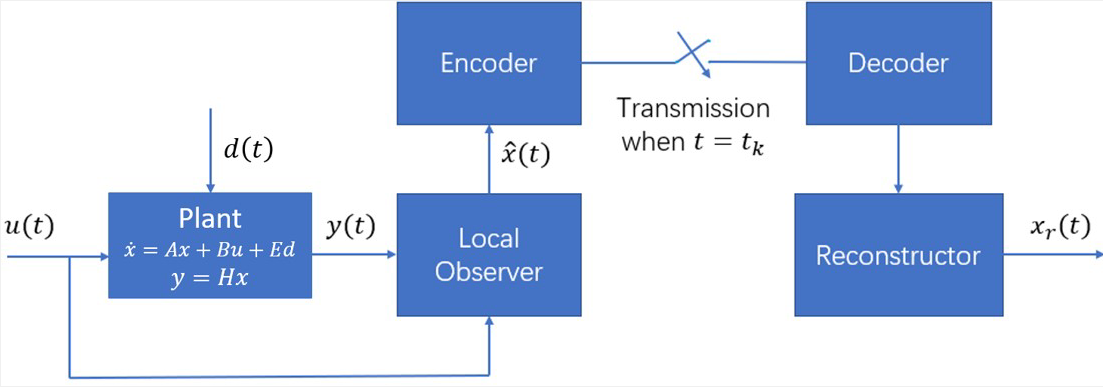}}
\caption{Remote state estimation setup}
\label{fig:detailed setup}
\end{figure}

We consider the problem of remote state estimation over a finite data rate channel in the setup depicted in Figure \ref{fig:detailed setup}. The plant has dynamics
\begin{equation}
    \label{dplant}
    \begin{split}
            \dot{x} &= Ax + Bu + Ed,\\
    y &= Hx,
    \end{split}
\end{equation}
with state $x \in \mathbb{R}^n$ input $u\in \mathbb{R}^m$,  output $y\in \mathbb{R}^{n_y}$, and unknown disturbance $d\in \mathbb{R}^{o}$. The system matrices $A, B, H, E$ are known, real matrices with appropriate dimensions. We assume that the pair $(A,H)$ is observable.

\begin{assumption} \label{ass:initial x u}
The initial state $x(0)$ and the input signals $Bu(t)$ and disturbance $Ed(t)$ for $t\in\mathbb{R}_{\geq 0}$ reside in known hypercubes $\mathcal{X} \subset \mathbb{R}^n$, $ \mathcal{U}\subset \mathbb{R}^n$ and $\mathcal{D}\subset \mathbb{R}^n$, respectively, defined as
\begin{equation}\label{initialstateset}
    \mathcal{X}:=  \mathcal{B}(x_c,x_b),\;
    \mathcal{U}:=  \mathcal{B}(\textbf{0}_m,u_b), \; \mathcal{D}:=\mathcal{B}(\textbf{0}_{o},d_b),
\end{equation}
where $x_c\in\mathbb{R}^n$ is the center of the initial set $\mathcal{X}$, and $x_b, u_b, d_b \in\mathbb{R}_{\geq 0}$.
\end{assumption}

Assumption \ref{ass:initial x u} requires the unknown initial state $x(0)$ of the system to reside in a known set $\mathcal{X}$ and the input $u$ to be bounded in $\mathcal{U}$, which is reasonable in most physical systems. 

\subsection{Local observer}
The local observer estimates the state of the plant \eqref{dplant} which is designed as follows 
\begin{equation}\label{observer dynamic}
    \dot {\hat{x}} = A\hat{x} +Bu + K(H\hat{x}-y),
\end{equation}
where $\hat{x} \in \mathbb{R}^n$ is the local state estimate. The local observer matrix $K\in\mathbb{R}^{n\times n_y}$ can be designed according to the following conditions.
\begin{assumption}
    \label{assum:observer}
There exist a matrix  $P = P^T\succ 0$, a matrix $Q$ and two scalars $\nu_1>0, \nu_2>0$ such that
\begin{equation}
\label{LMI2}
    \left [
    \begin{matrix}
    {A}^TP + P{A} + HQ^T + QH + \nu_1 I_{n} &P\\
       P               & -\nu_2 I_{n}
    \end{matrix}
    \right ] \preceq 0.
\end{equation}
\end{assumption}
 The inequality \eqref{LMI2} is a linear matrix inequality (LMI) in $P$, $Q$, $\nu_1$ and $\nu_2$ which can be solved efficiently with computational tools such as the LMI toolbox in MATLAB. The local observer matrix $K$ is then designed as $K = P^{-1}Q$. 

We assume that the local observer \eqref{observer dynamic} is initialized at $x_c$ (recall that $x_c \in \mathbb{R}^n$ comes from Assumption \ref{ass:initial x u}), i.e.,
\begin{equation}
    \label{local initial}
    \hat{x}(0) = x_c.
\end{equation}
Therefore, the initial local estimation error resides in the hypercube $\mathcal{X}$ from Assumption \ref{ass:initial x u}, i.e., $ \hat{e}(0) := x(0)-x_c \in \mathcal{X} $. This is essential as we use the bound on the local state estimation error $\hat{e}$ in approximating the reachable set of the local observer \eqref{local error dynamic} for the construction of the quantization region $\mathcal{S}_{Q}^{k}$ at each $t_k$, $k\in\mathbb{N}$. The proof of the following lemma can be found in Appendix \ref{proof:asp3}.

\begin{lemma}\label{lemma: asp3}
Consider the plant \eqref{dplant} and local observer \eqref{observer dynamic} under Assumptions \ref{ass:initial x u} - \ref{assum:observer}. The local state estimation error $\hat{e}(t) := x(t) - \hat{x}(t)$ of the local observer \eqref{observer dynamic} satisfies
\begin{equation}
    \label{ISS  error}
        |\hat{e}(t)| \leq \hat{\beta}(x_b,t) + \hat{\gamma}(d_b),\qquad  t \in \mathbb{R}_{\geq 0},
\end{equation}
where $x_c, x_b$ come from \eqref{initialstateset}, $\hat{\beta}(r,s) : = \sqrt{ \frac{n\lambda_{max}(P)}{\lambda_{min}(P)}}e^{-\frac{\lambda_es}{2}}r$  and $\hat{\gamma}(r): = \sqrt{ \frac{n\nu_2}{\lambda_{min}(P)\lambda_e}}r$, with $\lambda_e: = \frac{\nu_1}{n\lambda_{max}(P)}$; $\nu_1$ and $\nu_2$ come from \eqref{LMI2}. 
\end{lemma}

\begin{corollary}\label{coro lemma1}
    The set $\mathcal{E}_t$ of all the possible state estimation errors $\hat{e}(t)$ is given by
\begin{equation}\label{lemma asp3 1-1}
    \hat{e}(t) \in \mathcal{E}_t : = \mathcal{B}(\textbf{0}_n,\beta_d(t))
\end{equation}
for all $t \in \mathbb{R}_{\geq 0}$ with $\beta_d(t) := \hat{\beta}(x_b,t)+\hat{\gamma}(d_b)$, $\hat{\beta},
\hat{\gamma}$ are defined in \eqref{ISS  error}.
\end{corollary}

\subsection{Transmission times and the dynamic quantization scheme} \label{sec:quan_scheme}
We assume that the inter-transmission times are known and is periodic with period $T\in\mathbb{R}_{>0}$, i.e., the time at the $k$-th transmission satisfies $t_k=kT$, $k\in\mathbb{Z}_{+}$.

Given a quantization level\footnote{The quantization level $N$ is related to the number of bits $B_r\in\mathbb{Z}_{+}$ available to the communication channel according to $N =  {2^{B_r/n}}$, where $n$ is the dimension of the packet vector $x$. } $N\in\mathcal{Z}_{+}$, we employ a dynamic quantization scheme by updating the quantization region denoted by $\mathcal{S}^{k}_{Q} \subset \mathbb{R}^{n}$ at each transmission time $t_k$, $k\in\mathbb{Z}_{+}$. We perform reachability analysis to update the quantization region $\mathcal{S}^{k}_{Q}$. In Section \ref{sec:dynamic_set}, we propose using a set-based method of approximating the reachable set, which we achieve with zonotopes. We then compare it to a norm-based method of reachable set approximation in Section \ref{sec:dynamic_norm}. In both cases, we further regularize the reachable set as a hyperrectangle with $C^k\in\mathbb{R}^{n}$ being the center of the quantization region (the centroid), and $L^k\in\mathbb{R}^{n}$ being the segment length vector (the quantization range), i.e.,
\begin{equation}
    \label{Q reigon 0}
    \mathcal{S}^{k}_Q : = \mathcal{H}(C^k, L^k).
\end{equation}
We do so for ease of dividing $\mathcal{S}^{k}_Q$ into hyperrectangular sub-regions denoted by $\mathcal{S}_{q,j}^{k}$, $j\in\mathcal{I}_{n^N}$. Although other quantizer designs such as logarithmic or the more general Voronoi quantizer exist \cite{bullo2006quantized}, we consider the uniform quantizer in this paper for a fair comparison between the two reachable set approximation methods we will present later in Section \ref{sec:compare}. To this end, with quantization level $N\in\mathbb{Z}_{+}$, the quantization region $\mathcal{S}_{Q}^{k}$ is divided into $n^N\in\mathbb{Z}_{+}$ sub-regions at each $t_k$, $k\in\mathbb{Z}_{+}$. Each subregion is then numbered from $0$ to $N-1$ in the $i$-th dimension. Hence, the local state estimate $\hat{x}(t_k)$ is encoded into $\{0,1,\dots,N-1\}^{n}$, i.e, component $i\in\mathcal{I}_{n}$ of the encoded packet is defined by
\begin{equation} \label{eq:enco}
    P_{e,i}^{k} = Q_e(\hat{x}_{i}(t_k),C^{k}_{i},L^{k}_{i},N),
\end{equation}
where the encoding map $Q_e: \mathbb{R} \times\mathbb{R} \times\mathbb{R} \times\mathbb{Z}_+ \rightarrow \{0,1,...,N-1  \}$ is
\begin{equation}\label{eq:enco_map}
           Q_e(\hat{x}_i(t_k),C^{k}_{i},L^{k}_{i},N):= \left\lfloor{(\hat{x}_i(t_k)+L^{k}_i-C^{k}_i)\frac{N}{2L^{k}_i}}\right\rfloor.
\end{equation}
The parameters $L^{k}_{i},  C^{k}_{i}$ are the $i$-th element of the vector $L^{k}$ and $C^{k}$, respectively. 

Each subregion has a centroid $c_{j}^k \in \mathcal{S}_{Q}^{k} \subset \mathbb{R}^n$ with a quantization range $l^k:= \frac{L^k}{N}$, which is defined as
\begin{equation} \label{eq:subregion}
            \mathcal{S}_{q,j}^k := \mathcal{H}(c_{j}^k,l^k), \qquad j\in\mathcal{I}_{n^N}.
\end{equation} 
In our communication scheme, the decoded packet $P_d^{k}=(P_{d,1}^{k},P_{d,2}^{k},\dots,P_{d,n}^{k})\in\mathbb{R}^n$ is a member of the set of centroids $\{c_j^k\}_{j\in\mathcal{I}_{n^N}}$ of the subregions. Therefore, component $i \in \mathcal{I}_n$ of the decoded packet is given by 
\begin{equation} \label{eq:deco}
    P_{d,i}^{k} = Q_d(P^{k}_{e,i},C^{k}_i,L^{k}_i,N),
\end{equation} 
where the decoding map $Q_d: \mathbb{Z}_+ \times\mathbb{R} \times\mathbb{R} \times\mathbb{Z}_+ \rightarrow \mathbb{R}$ is 
\begin{equation}\label{eq:deco_map}
\begin{split}
           Q_d(P^{k}_{e,i},C^{k}_i,L^{k}_i,N) :=  C^{k}_{i}-\frac{L_i^k}{2}+\frac{L_i^k}{2N}(2P_{e,i}^{k}+1).
\end{split}
\end{equation}

\begin{remark}
We assume that the quantization parameters $C^k$ and $L^k$ are not transmitted over the communication channel. Instead, they are updated on both the encoder and the decoder sides simultaneously and therefore do not occupy any bandwidth. The robustness of a zoom-in and zoom-out dynamic quantization scheme in the context of stabilizing a linear system when the encoder and decoder asynchronously update their parameters was investigated in \cite{kameneva2009robustness}, but not in the context of the reachable-set based dynamic quantization scheme considered in this paper. 
\end{remark}

To prevent overflow at the first quantization step ($k = 0$), the initial quantization parameters satisfy the following.
\begin{assumption}\label{ass:initial_Q} 
The initial quantization region is chosen as $\mathcal{S}_Q^0 = \mathcal{H}(C^0,L^0)$ with $C^0 = x_c$, $L^0 = \textbf{1}_n x_b$, where $x_c$ and $x_b$ come from Assumption \ref{ass:initial x u}. 
\end{assumption}
Assumption \ref{ass:initial_Q} guarantees that the initial local state estimate $\hat{x}(0) \in \mathcal{S}_Q^0$ with arbitrary $l_0$. Moreover, we will show in Sections \ref{sec:dynamic_set} and \ref{sec:dynamic_norm} that by updating $C^k$ and $L^k$ for $k \in \mathbb{Z}_+$ appropriately, we can always prevent overflow, i.e., $\hat{x}(t) \in \mathcal{S}_Q^{k+1}$, for $t\in[t_k,t_{k+1})$, $k\in\mathbb{Z}_{+}$ as stated below.
\begin{assumption} \label{assum:contain}
    The quantization parameters $C^k$ and $L^k$ for $k \in \mathbb{Z}_+$ are updated such that $\hat{x}(t) \in \mathcal{S}_Q^{k+1}$ for $t\in[t_k,t_{k+1})$.  
\end{assumption}

\subsection{Reconstructor}
The reconstructor is responsible for producing a continuous prediction of the state $x$, based on the information received at discrete instances in time $t_k$. The dynamics of the reconstructed state $x_r \in \mathbb{R}^n$ satisfies, for all $k \in \mathbb{Z}_+$,
\begin{equation}\label{reconbehave}
\begin{split}
                \Dot{{x}}_r &= A{x}_r + K_r( Hx_r - HP_d^k  ),\ \forall t \in [t_k,t_{k+1}), \\
        {x}_r(t_{k}) &= P^{k}_d, 
\end{split} 
\end{equation}
with the reconstructor gain matrix $K_r\in\mathbb{R}^{n\times n_y}$ is chosen such that to satisfy $K_r=P^{-1}Q$, where matrices $P$ and $Q$ come from Assumption \ref{assum:observer}. The reconstructor \eqref{reconbehave} is initialized according to $x_r(0) =  x_c$, where $x_c \in \mathbb{R}^n$ comes from Assumption \ref{ass:initial x u}. At every transmission instant $t_k$, the reconstructor is initialised at $P_d^k$, which we recall is the decoded state at $t_k$.

\subsection{Objective and Approach}
The overall objective is to design a dynamic quantization scheme for the encoder and decoder, as well as a local observer and a reconstructor that converts the discrete-time decoder data into continuous-time reconstructed state $x_r(t) \in \mathbb{R}^n $, such that the state reconstruction error $e_r(t): = x(t) - x_r(t)$ is upper-bounded with respect to the quantization error $e_q(t):=\hat{x}(t)-P_d^k$, input $u$, and the unknown disturbance $d$. 

To this end, we provide the following guarantee for the state reconstruction error $e_r$ of our remote state estimation setup depicted in Figure \ref{fig:detailed setup}. Its proof can be found in Appendix \ref{proof:reconstruct_error}.
\begin{theorem} \label{thm:reconstruct_error}
    Consider the plant \eqref{dplant}, the local observer \eqref{observer dynamic}, encoder \eqref{eq:enco}, decoder \eqref{eq:deco} and reconstructor \eqref{reconbehave} under Assumptions \ref{ass:initial x u}-\ref{assum:contain}. The state reconstruction error $e_r(t): = x(t) - x_r(t)$ satisfies the following 
    \begin{align} \label{eq:recon_error_bound}
        |e_r(t)|\leq &\beta_r(|e_r(t_k)|,t-t_k) \notag \\
        &+\gamma_r\left(\max\{|u|_{[t_k,t]},|d|_{[t_k,t]},|\hat{e}|_{[t_k,t]},|{e}_q|_{[t_k,t]}\}\right),
    \end{align}
    for $t\in [t_k,t_{k+1})$ and for all $e_r(t_k)\in\mathbb{R}^n$, where $\beta_r\in\mathcal{KL}$, $\gamma_r\in\mathcal{K}$ and $|z|_{[t_k,t]}$ denotes $\sup_{s\in[t_k,t]}|z(s)|$.
\end{theorem}

From Theorem \ref{thm:reconstruct_error}, we see that the reconstruction error $e_r$ is ultimately bounded by the input $u$, disturbance $d$, local estimation error $\hat{e}$ and the quantization error $e_q$. Since we have established that $u$ and $d$ are bounded by Assumption \ref{ass:initial x u} and $\hat{e}$ by Lemma \ref{lemma: asp3}, we focus on establishing that the quantization error $e_q$ is also bounded. 

In this paper, we employ a reachable set-based approach in the dynamic quantization scheme. We propose over-approximating the reachable set with zonotopes in Section \ref{sec:dynamic_set} and compare it to a norm-based over-approximation of the reachable set in Section \ref{sec:dynamic_norm}. Both methods employ the idea of over-approximating the terminal reachable sets of the local state estimate $\hat{x}(t_k)$ and using it as the quantization region $S_{Q}^{k}$. To the best of our knowledge, the set-based approach using zonotopes which we will present in Section \ref{sec:dynamic_set}, has not been done in the literature. We show that this provides relaxed conditions over the norm-based dynamic quantization scheme in Section \ref{sec:compare}.

For both schemes, we compute the quantization error $e_q(t):=\hat{x}(t)-P^{k}_{d}$ and show that it is bounded for all $t\in\mathbb{R}_{\geq 0}$.        
\section{Set-based dynamic quantization using zonotopes}\label{sec:dynamic_set}
First, the dynamics of the local observer \eqref{observer dynamic} can be rewritten as
\begin{equation}
\label{observer dynamic 2}
    \dot {\hat{x}}  = A\hat{x} +Bu -KH \hat{e}.
\end{equation}
We can propagate the terminal reachable set of the local state estimate $\hat{x}$ based on the dynamics \eqref{observer dynamic 2} by treating the term $KH \hat{e}$ as an additional input. Since
\begin{equation*}
        KH \hat{e}(t )\in  KH\mathcal{E}_t \subseteq KH\mathcal{E}_{t_k},\ \forall t \in [t_k,t_{k+1}),
\end{equation*}
where $\mathcal{E}_t = \mathcal{B}(\textbf{0}_n,\beta_d(t))$ is defined in \eqref{lemma asp3 1-1}, and $\beta_d$ satisfies $\beta_d(t') \geq \beta_d(t'')$, for all $t'' \geq t' \geq 0$. Then 
\begin{equation*}
    Bu(t) - KH\hat{e}(t) \in \mathcal{U}\oplus KH\mathcal{E}_{t_{k}},\ \forall t\in [t_k,t_{k+1}).
\end{equation*}

Finally, given that $\hat{x}(t_k)\in \mathcal{S}_{Q}^{k}$, it implies that $\hat{x}(t_k)\in \mathcal{S}_{q,j}^{k}$ for some $j\in\mathcal{I}_{n^N}$. This is guaranteed for $k=0$ by Assumption \ref{ass:initial_Q}. Then, the terminal reachable set of the local state estimate $\hat{x}$ at the next transmission time $t_{k+1}$ can be deduced from Lemma \ref{2009}, that for $t\in[t_k,t_{k+1})$
\begin{equation}
    \label{eq:proga2}
        \hat{x}(t) \in \mathcal{R}_{[t_{k},t_{k+1}]}(\mathcal{S}_{q,j}^k, \mathcal{U}\oplus KH\mathcal{E}_{t_{k}})\subseteq \Lambda \mathcal{S}_{q,j}^k \oplus \mathcal{B}(\textbf{0}_n,\beta_{ue}^k) ,
\end{equation}
with $\Lambda:= e^{AT}$, where we recall that $T\in\mathbb{R}_{>0}$ is the inter-transmission interval, and $\beta_{ue}^k$ is the upper bound of the input term $(Bu - KH\hat{e})$ of the local observer dynamics \eqref{observer dynamic 2} for all $t \in [t_k,t_{k+1})$, which is
\begin{equation}\label{eq:betaue}
    \beta_{ue}^k: =|A|^{-1}e^{|A|T}(u_b + |KH|\beta_d(t_{k})),
\end{equation}
where $u_b$ comes from \eqref{initialstateset} and $\beta_d$ is defined in \eqref{ISS  error}.

We then over-approximate the terminal reachable set of the state estimate with a hyperrectangle such that the quantization region $\mathcal{S}^{k+1}_Q$ at the next transmission time $t_{k+1}$ is
\begin{equation}\label{eq:propa_2_end}
    \mathcal{S}_Q^{k+1} := \mathcal{H}(C^{k+1},L^{k+1})  \supseteq \left(  \Lambda \mathcal{S}_{q,j}^k \oplus \mathcal{B}(\textbf{0}_n,\beta_{ue}^k) \right).
\end{equation}
Hence, the dynamic quantization update law can be formulated as
\begin{equation}\label{eq:rangeupdate_zono}
    C^{k+1} = \Lambda P_d^{k},  \quad {L}_{i}^{k+1} = \beta_{ue}^k +\sum_{j=1}^n|\sigma_i^T\Lambda\sigma_j| \frac{{L}_j^k}{N}, \, i\in\mathcal{I}_{n},
\end{equation}
with initialization $P_d^0:=x_c$, $L_i^0:=x_b$ for $i\in \mathcal{I}_{n}$ according to Assumption \ref{ass:initial_Q} and we recall that $\beta_{ue}^k$ is defined in \eqref{eq:betaue}. We guarantee that by updating the dynamic quantization parameters $C^k$ and $L^k$ for $k\in\mathbb{Z}_{+}$ according to \eqref{eq:rangeupdate_zono}, the quantization region $S_{Q}^{k+1}$ always contains the local state estimate $\hat{x}(t)$ for $t\in[t_k,t_{k+1})$, $k\in\mathbb{Z}_{+}$, which we state in the following lemma.
\begin{lemma} \label{lem:dynamic_zono}
    Consider the plant \eqref{dplant}, the local observer \eqref{observer dynamic}, encoder \eqref{eq:enco}, decoder \eqref{eq:deco} under Assumptions \ref{ass:initial x u}-\ref{ass:initial_Q}. Suppose the quantization parameters $C^k$ and $L^k$ for $k\in\mathbb{Z}_{+}$ are updated according to the set-based scheme \eqref{eq:rangeupdate_zono}. Then Assumption \ref{assum:contain} holds. 
\end{lemma}
\begin{proof}
Let $k\in\mathbb{Z}_{+}$. Continuing from \eqref{eq:propa_2_end}, using standard properties of zonotopes and hyperrectangles outlined in Section \ref{sec:set}, we obtain 
\begin{align}
    \Lambda  S_{q,j}^{k} \oplus \mathcal{B}(\mathbf{0}_{n},\beta_{ue}^{k}) &  = \Lambda \mathcal{H}(c_j,l^k) \oplus \mathcal{B}(\mathbf{0}_{n},\beta_{ue}^{k}) \notag \\
    & = \Lambda \zono{c_j^k,G^k_{l}} \oplus \zono{\mathbf{0}_{n},G_{\beta}^{k}} \notag \\
    & = \zono{\Lambda c_j^k,(\Lambda G_l^{k}, G_{\beta}^{k} )},
\end{align}
where $G_{l}^{k}:=(\sigma_1^Tl^k\sigma_1,\sigma_2^Tl^k\sigma_2,\dots,\sigma_n^Tl^k\sigma_n)$ and $G_{\beta}^{k}:=(\beta_{ue}^k\sigma_1,\beta_{ue}^k\sigma_2,\dots,\beta_{ue}^k\sigma_n)$. Further, given \eqref{eq:propa_2_end}, we apply the containment property for zonotopes stated in Section \ref{sec:set} and obtain \eqref{eq:propa_2_end} where we note that $\sum_{m=1}^{n} |\sigma_i\Lambda\sigma_m^Tl^k\sigma_m|=\sum_{m=1}^{n} |\sigma_i\Lambda\sigma_m^T|l^k_m$ and $\sum_{m=1}^{n}|\sigma_i^T\beta_{ue}^{k}\sigma_m|=\beta_{ue}^{k}\sum_{m=1}^{n}|\sigma_i^T\sigma_m|=\beta_{ue}^k$ as $\sigma_i^T \sigma_m = 1$ only when $i=m$ and is $0$ otherwise. 
\end{proof}

We denote the quantization error at $t=t_k$ as $e_{q}^{k}:=\hat{x}(t_k)-P^{k}_{d}$. We show that its $i$-th component is always bounded by $\bar{e}^{k}_{q,i}$, i.e., 
\begin{equation} \label{eq:max_error_def}
    |e^{k}_{q,i}|\leq \bar{e}^{k}_{q,i}, \qquad  k\in\mathbb{Z}_{+}.  
\end{equation}
We call $\bar{e}^{k}_{q}:=(\bar{e}^{k}_{q,1},\bar{e}^{k}_{q,2}\dots, \bar{e}^{k}_{q,n})^T$ the maximum quantization error which we will show is bounded. Therefore, in conjuction with Lemma \ref{lem:dynamic_zono}, we guarantee that the quantization error $e_{q}(t)$ is also bounded in between transmission times, i.e. for $t\in[t_k,t_{k+1})$.

 By Lemma \ref{lem:dynamic_zono}, we have that $\hat{x}(t)\in \mathcal{S}_{Q}^{k+1}$ for $t\in[t_k,t_{k+1}]$, $k\in\mathbb{Z}_{+}$ which implies that $\hat{x}(t_k)\in \mathcal{S}_{q,j}^{k}$ for some $j\in\mathcal{I}_{n^N}$, where $\mathcal{S}_{q,j}^{k}$ is defined in \eqref{eq:subregion}. Hence, the maximum quantization error satisfies
\begin{equation} \label{eq:max_quan_error_def_l}
    \bar{e}_{q}^{k} = l^{k} = \frac{L^k}{N}, \qquad k\in\mathbb{Z}_{+}.
\end{equation}
Therefore, by the dynamic quantization law in \eqref{eq:rangeupdate_zono}, the maximum quantization error satisfies
\begin{equation} \label{eq:max_qerror_zono_dynamics}
    \bar{e}_{q}^{k+1} = \frac{\bar{\Lambda}}{N} \bar{e}_{q}^{k} + \frac{\beta_{ue}^{k}}{N} \mathbf{1}_{n},
\end{equation}
where 
\begin{equation} \label{eq:bar_lambda}
    \bar{\Lambda}:=(|\Lambda_{ij}|)_{i,j\in\mathcal{I}_{n}},
\end{equation}
is a matrix where each $(i,j)$-th component is the absolute value of $\Lambda_{ij}$ and $\Lambda_{ij}$ is the $(i,j)$-th component of matrix $\Lambda$ (recall from \eqref{eq:proga2} that $\Lambda=e^{AT}$).

We are now ready to show that by choosing the transmission time $T\in\mathbb{R}_{\geq 0}$ and total quantization level $N\in\mathbb{Z}_{>0}$ appropriately, the maximum quantization error $\bar{e}_{q}^{k}$ is bounded for all $k\in\mathbb{Z}_{+}$.
\begin{theorem}
     \label{prop:quan_error_zono}
Consider the plant \eqref{dplant}, the local observer \eqref{observer dynamic}, encoder \eqref{eq:enco}, decoder \eqref{eq:deco} under Assumptions \ref{ass:initial x u}-\ref{ass:initial_Q}. Suppose the quantization parameters $C^k$ and $L^k$ for $k\in\mathbb{Z}_{+}$ are updated according to \eqref{eq:rangeupdate_zono}. If the transmission interval $T>0$ and quantization level  $N\in\mathbb{Z}_{+}$ are chosen such that the matrix $\frac{\bar{\Lambda}}{N}$ is Schur, i.e.,
\begin{equation} \label{eq:cond_set}
    \max_{i\in\mathcal{I}_{n}}\left| \lambda_{i}\left(\frac{\bar{\Lambda}}{N}\right) \right| < 1,    
\end{equation}
 then the quantization error ${e}_{q}(t)$ is bounded for all $t\in[t_k,t_{k+1})$, $k\in\mathbb{Z}_{+}$.
\end{theorem}
\begin{proof}
    Consider the maximum quantization error system given in \eqref{eq:max_qerror_zono_dynamics}. Let  $T\in\mathbb{R}_{>0}$ and $N\in\mathbb{Z}_{>0}$ be chosen such that the matrix $\frac{\bar{\Lambda}}{N}$ is Schur. Therefore, the maximum quantization error system \eqref{eq:max_qerror_zono_dynamics} is asymptotically stable with $0$-input. We now show that the input $\beta_{ue}^{k}$ as defined in \eqref{eq:betaue} is bounded for all $k\in\mathbb{Z}_{+}$, i.e.,  $\beta_{ue}^{k} \leq \bar{\beta}_{ue}$ with $\bar{\beta}_{ue}:=|A|^{-1}(e^{|A|T}-1)(u_b + |KH|(\hat{\beta}(x_b,0)+\hat{\gamma}(d_b))$, where  $u_b$, $x_b$, $d_b \geq 0$ come from Assumption \ref{ass:initial x u}; $\hat{\beta}\in\mathcal{KL}$ and $\hat{\gamma}\in\mathcal{K}$ come from Lemma \ref{lemma: asp3}.  Therefore, the maximum quantization error $\bar{e}_{q}^k$ generated by system \eqref{eq:max_qerror_zono_dynamics} is bounded for all $k\in\mathbb{Z}_{+}$. Since $|e_{q,i}^k|\leq \bar{e}_{q,i}$, the quantization error $e_{q}^k$ is also bounded for all $k\in\mathbb{Z}_{+}$. Therefore, in combination with Lemma \ref{lem:dynamic_zono} where we have that $\hat{x}(t)\in S_{Q}^{k+1}$ for $t\in[t_{k+1},t_{k+1})$ $k\in\mathbb{Z}_{+}$, we can conclude that the quantization error $e_{q}(t)$ is bounded for $t\in[t_{k+1},t_{k+1})$, $k\in\mathbb{Z}_{+}$.
\end{proof}

\section{Norm-based dynamic quantization}\label{sec:dynamic_norm}
We now present a norm-based scheme in updating the quantization region $S_Q^{k}$, $k\in\mathbb{Z}_{+}$. We adapted the dynamic quantization scheme in \cite{de2004stabilizability} which was used in a stabilization setting to our remote state estimation problem. We choose the quantization region to be a hypercube, i.e., $S_{Q}^{k}:=\mathcal{B}(C(t_k),L^k)$, where the centroid $C$ has the following dynamics
\begin{align} \label{eq:norm_C_dynamics}
    \dot{C}(t) & = AC(t), \; t\in[t_k,t_{k+1}), \qquad C(t_k)  = P_d^k,
\end{align}
and we recall that $P_d^k$ is the decoded packet defined in \eqref{eq:deco}.

At the next transmission instant $t_{k+1}$, the quantization region $S_Q^{k+1}$ is designed as
\begin{equation} \label{eq:quan_region_norm}
    S_{Q}^{k+1} := \mathcal{B}(C(t_{k+1}),L^{k+1}),
\end{equation}
where the quantization parameters are updated according to
\begin{equation} \label{eq:rangeupdate_norm}
    C^{k+1}:=C(t_{k+1})=\Lambda P_d^{k}, \; L^{k+1} = \frac{e^{|A|T}}{N} L^k + \beta_{ue}^{k},    
\end{equation}
with initialization $P_d^0:=x_c$, $L^0:=x_b$, where $x_c$ and $x_b$ come from Assumption \ref{ass:initial x u}. Just as in the set-based method in Section \ref{sec:dynamic_set}, we show that the quantization region $S_{Q}^{k+1}$ contains the local state estimate $\hat{x}(t)$ for $t\in[t_k,t_{k+1})$, $k\in\mathbb{Z}_{+}$.

\begin{lemma} \label{lem:dynamic_norm}
    Consider the plant \eqref{dplant}, the local observer \eqref{observer dynamic}, encoder \eqref{eq:enco}, decoder \eqref{eq:deco}  under Assumptions \ref{ass:initial x u}-\ref{ass:initial_Q}. Suppose the quantization parameters $C^k$ and $L^k$ for $k\in\mathbb{Z}_{+}$ are updated according to the norm-based scheme \eqref{eq:quan_region_norm}-\eqref{eq:rangeupdate_norm}. Then Assumption \ref{assum:contain} holds.
\end{lemma}
\begin{proof}
    We update the quantization range $L^k$ based on the error $e_c:=\hat{x}-C$ between the local state estimate $\hat{x}$ and the centroid $C$ generated by \eqref{observer dynamic} and \eqref{eq:norm_C_dynamics}, respectively, which satisfies
\begin{align} \label{eq:ec_dynamics}
    \dot{e}_{c}(t) & = Ae_c(t) + Bu(t) - KH\hat{e}(t), \; t\in[t_k,t_{k+1}).
\end{align}
Therefore, $e_c$ can be bounded as follows for $t\in[t_k,t_{k+1})$,
\begin{align} \label{eq:e_c}
    |e_c(t)| & \leq e^{|A|T}|e_c(t_k)| \notag \\
    & \quad + \int_{t_k}^{t} e^{|A|(t-t_k-s)} \left( |Bu(s)| + |KH||\hat{e}(s)| \right) ds \notag \\
    & \leq e^{|A|T}|e_c(t_k)| + \beta_{ue}^{k},
\end{align}
where $\beta_{ue}^{k}$ is as defined in \eqref{eq:betaue} and we obtain the last inequality because $\int_{t_k}^{t} e^{|A|(t-t_k-s)} ds = |A|^{-1} e^{|A|T}\left( e^{-|A|t_k} - e^{-|A|t} \right) \leq |A|^{-1} e^{|A|T}\left( 1 - e^{-|A|t} \right) \leq |A|^{-1} e^{|A|T}$; $|Bu(s)| \leq u_b$ and $|\hat{e}(s)|\leq\beta_d(s)$ for $s\in\mathbb{R}_{\geq 0}$ according to Assumption \ref{ass:initial x u} and Lemma \ref{lemma: asp3}, respectively. 

We note that at $t=t_k$, we have that $|e_c(t_k)|:=|\hat{x}(t_k)-C(t_k)|=|\hat{x}(t_k)-P_d^k|$. Moreover, we have that  $|\hat{x}(t_k)-P_d^k|\leq \frac{L^k}{N}$. Thus, from \eqref{eq:e_c}, 
\begin{equation} \label{eq:error_bound_norm}
    |e_c(t)| \leq \frac{e^{|A|T}}{N} L^k + \beta_{ue}^{k}, \qquad t\in[t_k,t_{k+1}),\, k\in\mathbb{Z}_{+}.
\end{equation}

Therefore, we conclude that $$\hat{x}(t)\in\mathcal{B}\left(P_d(t_k),\frac{e^{|A|T}}{N} L^k + \beta_{ue}^{k}\right), \qquad t\in[t_k,t_{k+1}),$$ which satisfies \eqref{eq:quan_region_norm} and \eqref{eq:rangeupdate_norm}. 
\end{proof}

 We can now provide conditions on the transmission time $T\in\mathbb{R}_{\geq 0}$ and the total quantization level $N\in\mathbb{Z}_{>0}$ such that quantization error ${e}_{q}(t):=\hat{x}(t)-P_d^k$  is bounded for $t\in[t_k,t_{k+1})$. 
\begin{theorem} \label{prop:quan_error_norm}
    Consider the plant \eqref{dplant}, the local observer \eqref{observer dynamic}, encoder \eqref{eq:enco}, decoder \eqref{eq:deco}  under Assumptions \ref{ass:initial x u}-\ref{ass:initial_Q}. Suppose the quantization parameters $C^k$ and $L^k$ for $k\in\mathbb{Z}_{+}$ are updated according to the norm-based scheme \eqref{eq:quan_region_norm}-\eqref{eq:rangeupdate_norm}. If the inter-transmission interval $T\in\mathbb{R}_{\geq 0}$ and total quantization level  $N\in\mathbb{Z}_{>0}$ are chosen such that \begin{equation} \label{eq:cond_norm}
        \frac{e^{|A|T}}{N}<1,
    \end{equation} 
    then the quantization error ${e}_{q}(t)$ is bounded for all $t\in\mathbb{R}_{\geq 0}$.
\end{theorem}
\begin{proof}
    Recall that the quantization error at $t_k$ is $e_{q}^{k}:=\hat{x}(t_k)-P_d(t_k)$, which according to Lemma \ref{lem:dynamic_norm}, satisfies $|e_q^{k}| \leq L^k$ since the quantization region $S_Q^{k}$ is defined by \eqref{eq:quan_region_norm}. Therefore, we need to show that $L^k$ is bounded for all $k\in\mathbb{Z}_{+}$ to conclude boundedness of the quantization error $e_q^{k}$ for all $k\in\mathbb{Z}_{+}$. According to the update scheme for the quantization parameters \eqref{eq:rangeupdate_norm}, the range  $L_k$ satisfies
\begin{align} \label{eq:Lk_bound}
    L^k & \leq \frac{|A|T}{N} L^0 + \underset{m\in\mathbb{N}_{[0,k]}}{\sum} \left(\frac{|A|T}{N} \right)^{k-1-m} {\beta_{ue}^{0}}, 
\end{align}
where $\beta_{ue}^{0}$ satisfies $\beta_{ue}^{k}\leq \beta_{ue}^{0}$ for all $k\in\mathbb{Z}_{+}$. Recall that ${\beta_{ue}^{k}}\in\mathcal{L}$ is defined in \eqref{eq:betaue}. Therefore, the right-hand side of inequality \eqref{eq:Lk_bound} is only bounded when $\frac{|A|T}{N}<1$, which concludes the proof. 
\end{proof}

\subsection{Comparison} \label{sec:compare}
The difference between the set-based and norm-based schemes lies in how they over-approximate the `size' of the reachable set, which is captured by the quantization range $L^k$. The set-based method propagates the zonotope, and thereby provides upper bounds for each individual component in the state, while the norm-based method only considers the norm of $e^{AT}$. In fact, the set-based method provides a relaxed condition \eqref{eq:cond_set} in Theorem \ref{prop:quan_error_zono} compared to the norm-based method. We see this by first establishing the following. 

\begin{lemma} \label{lem:compare}
Given any matrix $A\in\mathbb{R}^{n\times n}$, scalars $T\in\mathbb{R}_{\geq 0}$ and  $N\in\mathbb{Z}_{>0}$, the following holds
\begin{equation} \label{eq:compare}
    \max_{i\in\mathcal{I}_{n}}\left|\lambda_i\left( \frac{\bar{\Lambda}}{N} \right)\right|  \leq \frac{e^{|A|T}}{N}, 
\end{equation}
where $\bar{\Lambda}$ is as defined in \eqref{eq:bar_lambda}.
\end{lemma}
\begin{proof}
    Using Theorem 5.6.9 in \cite{horn2012matrix}, we have that
    \begin{align} \label{eq:compare_int}
        \max_{i\in\mathcal{I}_{n}}\left|\lambda_i\left( \frac{\bar{\Lambda}}{N} \right)\right|  \leq \left| \frac{\bar{\Lambda}}{N}\right| =\frac{1}{N}|\Lambda|,
    \end{align}
    where we obtain the equality since $N\in\mathbb{Z}_{>0}$ and that by definition of the infinity norm of a matrix and the matrix $\bar{\Lambda}$ in \eqref{eq:bar_lambda}, respectively, we obtain
    \begin{align*}
        |\bar{\Lambda}| = \max_{j\in\mathcal{I}_{n}} \sum_{i\in\mathcal{I}_{n}}|\bar{\Lambda}_{ij}| = \max_{j\in\mathcal{I}_{n}} \sum_{i\in\mathcal{I}_{n}}|{\Lambda}_{ij}| = |\Lambda|.
    \end{align*}
    Therefore, from \eqref{eq:compare_int} and $|\Lambda| \leq e^{|A|T}$, we obtain \eqref{eq:compare}. 
\end{proof}

By Lemma \ref{lem:compare}, we see that condition \eqref{eq:cond_set} of the set-based scheme in Theorem \ref{prop:quan_error_zono} relaxes condition \eqref{eq:cond_norm} of the norm-based scheme in Theorem \ref{prop:quan_error_norm}. We state this result in the proposition below.

\begin{proposition} \label{prop:compare}
    Consider the plant \eqref{dplant}, the transmission interval $T\in\mathbb{R}_{\geq 0}$ and quantization level $N\in\mathbb{Z}_{>0}$. If the parameters $T$ and $N$ satisfy condition \eqref{eq:cond_norm} of the norm-based scheme, then condition \eqref{eq:cond_set} of the set-based scheme will also hold.
\end{proposition}

Consequently, when designing a dynamic quantization scheme, the designer first checks if the norm-based method (Theorem \ref{prop:quan_error_norm}) is met. If so, then the designer has the option of both the set and norm-based methods presented so far in Sections \ref{sec:dynamic_set} and \ref{sec:dynamic_norm}, respectively. Future work would include deriving analytical bounds on the quantization error for both methods. In this paper, we compare them in simulation which we present in the next section.

\section{Numerical simulation} \label{sec:sim}
We consider a $2$-dimensional LTI system in the form of \eqref{dplant} with 
\begin{equation}\label{simulationplant}
\begin{split}
        A = \begin{bmatrix}
    -1  & -4   \\
    4  & -1  
    \end{bmatrix},\ B =  E = 
    \begin{bmatrix}
    1\\
    1
    \end{bmatrix},\ H =
    \begin{bmatrix}
    1 & 0\\
    \end{bmatrix},
\end{split}
\end{equation}
with the initial state $x(0) \in \mathcal{X} =  \mathcal{B}([10,-5]^T,1) = [9,11]\times[-6,-4]$, 
the input $u(t) = 0.5\sin(t) \in \mathcal{B}(0,0.5),\ \forall t\geq 0$, and the disturbance $d(t)$ is a random noise that is upper-bounded by $0.05$ for all $t\geq 0$.
Hence, Assumption \ref{ass:initial x u}
 is satisfied with $x_c =[10,5]^T ,x_b = 1, u_b = 0.5$ and $d_b = 0.05$.

We solve \eqref{LMI2} using the LMI toolbox in {MATLAB} to obtain
\begin{equation}
\begin{split}
        P = \begin{bmatrix}
        2.0648 &0.9237\\
        0.9237 & 1.9195
    \end{bmatrix}&,\ Q = \begin{bmatrix}
        -7.7353 \\
        -0.0248
    \end{bmatrix},\\
    \nu_1 = 8.2561&,\ \nu_2 = 7.2571.
\end{split}
\end{equation}
Hence, the local observer gain $K$ of \eqref{observer dynamic} and the reconstructor gain $K_r$ of \eqref{reconbehave} are 
\begin{equation}
    K =K_r= P^{-1}Q = \begin{bmatrix}
        -4.7666\\
        2.2808
    \end{bmatrix}.
\end{equation}

In our simulation study, the plant \eqref{dplant}, the local observer \eqref{observer dynamic} and reconstructor \eqref{reconbehave} are initialized at $x(0) = [10.5, -5.5]^T$, $\hat{x}(0)  = x_c =[10,5]^T$ and ${x}_r(0)  = x_c =[10,5]^T$, respectively. 

Our communication setup depicted in Figure \ref{fig:detailed setup} has quantization level $N = 4$, bit rate $B_r = 4$, and the inter-transmission interval is $T = 0.1$. With these parameters, we can verify that $\frac{\bar{\Lambda}}{N}$ is Schur and $\frac{e^{|A|T}}{N}<1$, where each satisfies the condition of the set-based  (Theorem \ref{prop:quan_error_zono}) and norm-based (Theorem \ref{prop:quan_error_norm}) quantization schemes, respectively. Therefore, the set-based and norm-based dynamic quantization schemes presented in Sections \ref{sec:dynamic_set} and \ref{sec:dynamic_norm}, respectively, are applicable. Moreover, the reconstructor $x_r$ given by \eqref{reconbehave} has bounded reconstruction error according to Theorem \ref{thm:reconstruct_error}. We simulated for the time-interval $[0,20]s$, i.e., with transmission interval $T=0.1$, we transmit $200$ times. By defining the 
steady state quantization error as ${e}_{q}^{\infty}:=\lim_{k \to \infty } |{e}^k_{q}|$ and
steady-state reconstruction error as $\bar{e}_{r}^{\infty}:=\lim_{t\to \infty} |{e}_{r}(t)|$, we compare them for each scheme, which we summarise  in Table \ref{tab:compare}. In both cases, the set-based scheme outperforms the norm-based scheme.
\begin{table}[ht]
	\centering
	\caption{Set-based and norm-based dynamic quantization schemes}
	\label{tab:compare}
	\begin{tabular}{|p{0.18\textwidth}|c|c|}
		\hline 
		&\multicolumn{2}{|c|}{Dynamic quantization schemes}  \\		 
	   \cline{2-3}	
            &Set-based & Norm-based \\ 
		\hline 
		Steady-state quantization error ${e}^{\infty}_{q}$ & $0.0571$ & $0.0684$ \\ 
		\hline 
    Steady-state reconstruction error ${e}^{\infty}_{r}$ & 
    0.0921& 0.1170 \\
    \hline
	\end{tabular} 
\end{table}

\section{Conclusions} \label{sec:conclude}
We presented a reachability analysis approach to design the quantization region in a dynamic quantization scheme using a uniform quantizer, which is employed in a remote state estimation setting. Our results guarantee that the reconstruction error is upper bounded by a class $\mathcal{K}$ function of the quantization error. To this end, two different methods were presented: set-based and norm-based. The set-based method using zonotopes is shown to yield a less conservative sufficient condition over the norm-based method, for guaranteeing a bounded quantization error. Moreover, we also obtain in simulations that the set-based method is an improvement over the norm-based method, whereby its steady-state quantization error is smaller. 
\appendix

\subsection{Proof of Lemma \ref{lemma: asp3}} \label{proof:asp3}
The dynamics of the local estimation error $\hat{e}:={x} - \hat{x}$ satisfies, for all $t \in \mathbb{R}_+$,
\begin{equation}
    \label{local error dynamic}
    \dot{\hat{e}} = (A+K H)\hat{e} + Ed.
\end{equation}

First, we show that the observer \eqref{observer dynamic} satisfies \eqref{ISS  error} by choosing the observer matrix $K = P^{-1}Q$ and using a candidate Lyapunov function $V(\hat{e}) = \hat{e}^T P \hat{e}$, where $P = P^T\succ 0$ and matrix $Q$ satisfies \eqref{LMI2}. The time derivative of $V(\hat{e})$ along the trajectories of the state estimation error system \eqref{local error dynamic} is
\begin{equation*}\label{prof lemma2 1-1}
\dot{V}(\hat{e}) = \hat{E}^T
\left [
    \begin{matrix}
    (A+KH)^TP + P(A+KH)   &P\\
       P               & 0I_n
    \end{matrix}
    \right ]\hat{E},
    \end{equation*}
where $\hat{E}: = [\hat{e}^T, Ed^T]^T$. Since $K=P^{-1}Q$ and due to \eqref{LMI2}, we can further deduce that
\begin{equation}\label{vdot}
\dot{V}(\hat{e})  \leq -\nu_1|\hat{e}|^2 + \nu_2n|Ed|^2.
\end{equation}
Further, $V(\hat{e})$ can be bounded as follows 
\begin{equation}
    \label{sandwich2}
            \lambda_{min}(P)|\hat{e}|^2 \leq V(\hat{e}) \leq n\lambda_{max}(P)|\hat{e}|^2.
\end{equation}
By applying the right inequality of (\ref{sandwich2}), we can obtain that $\dot{V}(\hat{e})$ satisfies
\begin{equation}
\label{leq}
 \dot{V}(\hat{e}) \leq -\frac{\nu_1}{n\lambda_{max}(P)} V(\hat{e}) + \nu_2n|Ed|^2.
\end{equation}
Next, let $\lambda_e := \frac{\nu_1}{n\lambda_{max}(P)}$ and by the Comparison Principle, we have
\begin{equation}
\begin{split}
    {V}(\hat{e}(t)) 
        &\leq e^{-\lambda_e t}V(\hat{e}(0)) + \int_{0}^t e^{-\lambda_e (t-s)}\nu_2 n|Ed(s)|^2ds\\
        & \leq e^{-\lambda_e t}n\lambda_{max}(P)|\hat{e}(0)|^2 + \frac{\nu_2nd_b^2}{\lambda_e},
\end{split}
\end{equation}
where we obtain the last inequality because $|Ed(s)| \leq d_b$ for all $s\in \mathbb{R}_+$, $d_b$ according to Assumption \ref{ass:initial x u}, and $\int_{0}^t e^{-\lambda_e (t-s)}ds = (1-e^{-\lambda_e t})/\lambda_e \leq 1/\lambda_e$. Finally, by applying the left inequality of (\ref{sandwich2}), we obtain
\begin{equation} 
     |\hat{e}(t)|^2 \leq \frac{V(\hat{e}(t))}{\lambda_{min}(P)} 
     =  \frac{n\lambda_{max}(P)}{\lambda_{min}(P)}e^{-\lambda_e t}|\hat{e}(0)|^2 + \frac{\nu_2nd_b^2}{\lambda_{min}(P)\lambda_e}. \nonumber
\end{equation}
Since $ \sqrt{a + b} \leq \sqrt{a} + \sqrt{b}$ for all $a,b \in \mathbb{R}_{\geq 0}$, we get the following for $t\in \mathbb{R}_{\geq 0}$
\begin{equation}
\label{sqrtbound}
            |\hat{e}(t)| \leq \sqrt{ \frac{n\lambda_{max}(P)}{\lambda_{min}(P)}}e^{-\frac{\lambda_et}{2}}|\hat{e}(0)| + \sqrt{ \frac{\nu_2n}{\lambda_{min}(P)\lambda_e}}d_b.
\end{equation}
Since the local state estimate $\hat{x}(0)$ is initialized at $x_c$ according to \eqref{local initial}, it can be deduced that $|\hat{e}(0)|= |x(0) - \hat{x}(0)|  =  |x(0) - x_c|\leq x_b$. Therefore, we obtain \eqref{ISS  error}. \hfill $\Box$

\subsection{Proof of Theorem \ref{thm:reconstruct_error}} \label{proof:reconstruct_error}
The reconstruction error $e_r:=x-x_r$ system from \eqref{dplant} and \eqref{reconbehave} is
\begin{align}
    \dot{e}_r & = (A+K_rH)e_r + Bu + Ed -K_rH(\hat{e}+e_q),
\end{align}
where $\hat{e}:=x-\hat{x}$ and $e_q:=\hat{x}-P_d^k$. Since $K_r$ is chosen in the same manner as $K$ from the local observer \eqref{local error dynamic}, the proof follows closely the proof for Lemma \ref{lemma: asp3} using a quadratic Lyapunov function $V(e_r):=e_r^T P e_r$, where $P=P^T\succ 0$ satisfies Assumption \ref{assum:observer}. Thus, we conclude the proof by obtaining \eqref{eq:recon_error_bound} with $\beta_r(r,s):=\sqrt{\frac{n\lambda_{\max}(P)}{\lambda_{\min}(P)}}e^{\frac{\lambda_e}{2}s}r$, $\lambda_e:=\frac{\nu_1}{n \lambda_{\max}(P)}$, ${\gamma}_r(s):=\sqrt{\frac{\nu_2 n}{\lambda_{\min}(P)\lambda_e}}(|E|+|B|+|K_r H|) s$. 
\hfill $\Box$

\addtolength{\textheight}{-12cm}   

\bibliographystyle{ieeetr}

\end{document}